\documentclass[12pt]{amsart}
\usepackage{amssymb}
\theoremstyle{plain}
\newtheorem{thm}{Theorem}
\theoremstyle{definition}
\newtheorem{prop}{Proposition}[section]

\newtheorem{lem}[prop]{Lemma}

\newtheorem{defn}[prop]{Definition}
\newtheorem{claim}[prop]{Claim}
\theoremstyle{remark}
\newtheorem{subclaim}{Subclaim}

\DeclareMathOperator{\env}{env}
\newcommand{\model}{H(\chi)}
\DeclareMathOperator{\str}{Str} \DeclareMathOperator{\gm}{Gm}

\DeclareMathOperator{\bd}{bd} \DeclareMathOperator{\Sk}{Sk}

\DeclareMathOperator{\tcf}{tcf} \DeclareMathOperator{\otp}{otp}
\DeclareMathOperator{\acc}{acc} \DeclareMathOperator{\pr}{Pr}
 \DeclareMathOperator{\ran}{ran}
 
 \DeclareMathOperator{\cf}{cf}

\newcommand{\sk}{\vskip.05in}
\DeclareMathOperator{\id}{id}

\newcommand{\restr}{\upharpoonright}

\DeclareMathOperator{\nacc}{nacc}

\newcommand{\subs}{\subseteq}
\newcommand{\sups}{\supseteq}
 
\newcommand{\jbk}{J^{\bd}_\kappa}

\numberwithin{equation}{section}
\begin{document}
\title[Coloring Theorems]{Successors of singular cardinals and coloring theorems {I}}
\author{Todd Eisworth}
\address{Department of Mathematics\\
         University of Northern Iowa\\
         Cedar Falls, IA\\
         50614}
\email{eisworth@uni.edu}
\author{Saharon Shelah}
\address{Institute of Mathematics\\
    Hebrew University of Jerusalem\\
    Jerusalem, Israel}
\email{shlhetal@math.huji.ac.il}  \subjclass{}
\thanks{This is publication number 535 of the second author.}
\date{\today}
\begin{abstract}
We investigate the existence of strong colorings on successors of
singular cardinals. This work continues Section 2 of
\cite{Sh:413}, but now our emphasis is on finding colorings of
pairs of ordinals, rather than colorings of finite sets of
ordinals.
\end{abstract}
\maketitle
\section{Introduction}

The theme of this paper is that strong coloring theorems hold at
successors of singular cardinals of uncountable cofinality, except
possibly in the case where the singular cardinal is a limit of
regular cardinals that are Jonsson in a strong sense.

Our general framework is that $\lambda=\mu^+$, where $\mu$ is
singular of uncountable cofinality. We will be searching for
colorings of pairs of ordinals $<\lambda$ that exhibit quite
complicated behaviour.  The following definition (taken from
\cite{Sh:g}) explains what ``complicated'' means in the previous
sentence.

\begin{defn}
Let $\lambda$ be an infinite cardinal, and suppose
$\kappa+\theta\leq\mu\leq\lambda$.
$\pr_1(\lambda,\mu,\kappa,\theta)$ means that there is a symmetric
two--place function $c$ from $\lambda$ to $\kappa$ such that if
$\xi<\theta$ and for $i<\mu$, $\langle
\alpha_{i,\zeta}:\zeta<\xi\rangle$ is a strictly increasing
sequence of ordinals $<\lambda$ with all $\alpha_{i,\zeta}$'s
distinct, then for every $\gamma<\kappa$ there are $i<j<\mu$ such
that
\begin{equation}
\zeta_1<\xi\text{ and }\zeta_2<\xi\Longrightarrow
c(\alpha_{i,\zeta_1},\alpha_{i,\zeta_2})=\gamma.
\end{equation}
\end{defn}

Just as in \cite{Sh:413}, one of our main tools is a game that
measures how ``Jonsson'' a given cardinal is.

Recall that a cardinal $\lambda$ is a Jonsson cardinal if for
every $c:[\lambda]^{<\omega}\rightarrow\lambda$, we can find a
subset $I\subs\lambda$  of cardinality $\lambda$ such that the
range of $c\restr I$ is a proper subset of $\lambda$.  A reader
seeking more background should investigate \cite{Sh:365} and
\cite{Sh:380} in \cite{Sh:f}.

\begin{defn}
Assume $\mu\leq\lambda$ are cardinals, $\gamma$ is an ordinal,
$n\leq\omega$, and $J$ is an ideal on $\lambda$. We define the
game $\gm^n_J[\lambda,\mu,\gamma]$ as follows: \sk \noindent A
play lasts $\gamma$ moves. \sk \noindent In the
$\alpha^{\text{th}}$ move, the first player chooses a function
$F_\alpha:[\lambda]^{<n}\rightarrow\mu$, and the second player
responds by choosing (if possible) a subset $A_\alpha\subs\lambda$
such that
\begin{itemize}
\item $A_\alpha\subs\bigcap_{\beta<\alpha}A_\beta$
\sk
\item $A_\alpha\in J^+$
\sk
\item $\ran(F_\alpha\restr[A_\alpha]^{<n})$ is a proper subset of $\mu$.
\end{itemize}
\sk \noindent The second player loses if he has no legal move for
some $\alpha<\gamma$, and he wins otherwise.
\end{defn}

In the previous definition, if $J=J^{\bd}_\lambda$ then we may
omit it. Note that it causes no harm if we use a set $E$ of
cardinality $\lambda$ instead of $\lambda$ itself; in this case,
we write $\gm^n_J[E,\mu,\gamma]$.

Note that $\lambda$ is a Jonsson cardinal if and only if Player I
does not have a winning strategy in the game
$\gm^\omega[\lambda,\lambda,1]$.  One may view the lack of a
winning strategy for Player I in games of longer length as a
strong version of Jonsson-ness or a weak version of measurability
--- if $\lambda$ is measurable, then Player II can make sure her
moves are elements of some $\lambda$--complete ultrafilter.

The following claim investigates how the existence of winning
strategies is affected by modifications to the game; the proof is
left to the reader.
\begin{claim}
\label{1.4} \hfill
\begin{enumerate}
\item If $\mu'\leq\mu$ and the first player has a winning strategy in $\gm^n_J[\lambda,\mu,\gamma]$,
 then she has a winning strategy in $\gm^n_J[\lambda,\mu',\gamma]$.
\sk
\item Suppose we weaken the demand on the second player to
\begin{equation}
\text{``}(\exists\zeta<\lambda)[\ran(F_\alpha\restr[A_\alpha\setminus\zeta]^{<n})\text{
is a proper subset of $\mu$}].\text{''}
\end{equation}
If $\cf(\lambda)\geq\gamma$ and $J\sups J^{\bd}_\lambda$, then the
first player has a winning strategy in the revised game if and
only if she has a winning strategy in the original game. \sk
\item If $J$ is $\gamma$--complete, then the same applies to the case where we weaken the
 demand on the second player to
\begin{equation}
\text{``}(\exists Y\in J)[\ran(F_\alpha\restr[A_\alpha\setminus
Y]^{<n})\text{ is a proper subset of $\mu$}].\text{''}
\end{equation}
\sk
\item We can allow the second player to pass, i.e., to let $A_\alpha=\bigcap_{\beta<\alpha}A_\beta$
 (even if this is not a legal move) as long as we declare that the second player loses if the
  order--type of the set of moves where he did not pass is $<\gamma$.
\sk
\item \label{usedgame} If Player I has a winning strategy in $\gm^n_J[\lambda,\mu,\gamma]$ for
 every $\mu<\mu^*$ where $\mu^*$ is singular and $\gamma>\cf(\mu^*)$ is regular, then Player I has a winning strategy  in $\gm^n_J[\lambda,\mu^*,\gamma]$.
 We can weaken the requirement that $\gamma$ is regular and instead require that $\cf(\gamma)>\cf(\mu^*)$ and $\omega^\gamma=\gamma$.
\sk
\end{enumerate}
\end{claim}

In Section 2 of \cite{Sh:413}, the existence of winning strategies
for Player I in variants of the game is investigated. We will
prove one such result here; the reader should look in
\cite{Sh:413} for others.

\begin{claim}
\label{1.4b} If $2^\chi<\lambda<\beth_{(2^\chi)^+}(\chi)$ then
Player I has a winning strategy in
$\gm^\omega[\lambda,\chi,(2^\chi)^+]$.
\end{claim}
\begin{proof}
At a stage $i$, Player I will select a function
$F_i:[\lambda]^{<\omega}\rightarrow\chi$ coding the Skolem
functions of some model $M_i$.

For the initial move, we let the model $M_0$ have universe
$\lambda$, and include in our language all relations on $\lambda$
and all functions from $\lambda$ to $\lambda$ of any finite arity
that are first order definable in the structure $\langle
H(\lambda^+),\in,<^*_{\lambda^+}\rangle$ with the parameters
$\chi$ and $\lambda$.

For subsequent moves, $M_i$ is an expansion of $M_0$ with universe
$\lambda$ that has all relations on $\lambda$ and all functions
from $\lambda$ to $\lambda$ of any finite arity that are first
order definable in the structure $\langle
H(\lambda^+),\in,<^*_{\lambda^+}\rangle$ from the parameters
$\chi$, $\lambda$, $M_0$, and $\langle A_j:j<i\rangle$.

To obtain the function $F_i$, we let $\langle
F_n^i:n<\omega\rangle$ list the Skolem functions of $M_i$ in such
a way that $F_n^i$ has $m_i(n)\leq n$ places. Let
$h:\omega\rightarrow\omega$ be such that for all $n$, $h(n)\leq n$
and $h^{-1}(\{n\})$ is infinite.   We then define
\begin{equation}
F_i(u)=
\begin{cases}
F^i_{h(|u|)}(\{\alpha\in u:|u\cap\alpha|<m_i(n)\}) &\text{if this is $<\chi$}\\
0 &\text{otherwise}
\end{cases}
\end{equation}
The point of doing this is that whenever Player II chooses $A_i$,
we know that $\ran(F_i\restr[A_i]^{<\omega})$ will look like the
result of intersecting an elementary submodel of $M_i$ with
$\chi$; in particular, this range will be closed under the
functions from $M_i$.

Note that $M_0$ (and all expansions of it) has definable Skolem
functions and so for any $i$ and $A\subs\lambda$, the Skolem hull
of $A$ in $M_i$ (denoted by $\Sk_{M_i}(A)$) is well--defined.

Let $\langle (F_i, A_i):i<(2^\chi)^+\rangle$ be a play of the game
in which Player I uses this strategy (with $M_i$ the model
corresponding to $F_i$). For each $i$, define
\begin{equation}
\alpha_i=\min\{\alpha:\left|\Sk_{M_0}(A_i)\cap\beth_\alpha(\chi)\right|>\chi\}.
\end{equation}

By the choice of $M_0$ and $M_i$, clearly $\alpha(i)$ is a
successor ordinal or a limit ordinal of cofinality $\chi^+$, and
\begin{equation}
\label{assumption} |\Sk_{M_0}(A_i)\cap\beth_{\alpha_i}(\chi)|\leq
2^\chi.
\end{equation}

Since $A_i\subs A_j$ for $i>j$, we know the sequence $\langle
\alpha_i:i<(2^\chi)^+\rangle$ is non--decreasing. Furthermore, for
each $i$ we know
\begin{equation}
\alpha_i<\min\{\beta:\lambda\leq\beth_\beta(\chi)\}<(2^\chi)^+.
\end{equation}
This means that  the sequence
$\langle\alpha_i:i<(2^\chi)^+\rangle$ is eventually constant, say
with value $\alpha^*$.  Let $i^*$ be the least ordinal
$<(2^{\chi})^+$ such that $\alpha_i=\alpha^*$ for $i\geq i^*$.

\begin{prop}
If $i^*\leq i<(2^{\chi})^+$, then
$\Sk_{M_0}(A_{i+1})\cap\beth_{\alpha^*}(\chi)$ is a proper subset
of $\Sk_{M_0}(A_i)\cap\beth_{\alpha^*}(\chi)$.
\end{prop}

\begin{proof}
Note that $i^*$, $\alpha^*$, and $\beth_{\alpha^*}(\chi)$ are all
elements of $M_{i+1}$ as they are definable in $\langle
H(\lambda^+),\in,<_{\lambda^+}\rangle$ from the parameters $M_0$
and $\langle A_j:j\leq i\rangle$.  Furthermore,
\begin{equation}
\gamma^*:=\min\{\gamma<\lambda:|\Sk_{M_0}(A_i)\cap\gamma|=\chi\}
\end{equation}
is also definable in $M_{i+1}$ (and $<(2^\chi)^+$). Thus the
language of $M_{i+1}$ includes a bijection between
$\Sk_{M_0}(A_i)\cap\gamma^*$ and $\chi$.

If Player I has not won the game at this stage, after Player I
selects $A_{i+1}$ we will be able to find an ordinal $\beta<\chi$
such that $\beta\notin\ran(F_{i+1}\restr[A_{i+1}]^{<\omega})$.  By
definition of $h$, we know $\beta':=h^{-1}(\beta)$ is an element
of $\Sk_{M_0}(A_i)\cap \beth_{\alpha^*}(\chi)$.  However, $\beta'$
is not an element of $\Sk_{M_{i+1}}(A_{i+1})$ -- since $F_{i+1}$
codes the Skolem functions of $M_{i+1}$, the range of
$F_{i+1}\restr[A_{i+1}]^{<\omega}$ is
$\Sk_{M_{i+1}}(A_{i+1})\cap\chi$. Since $\Sk_{M_{i+1}}(A_{i+1})$
is closed under $h$, this  contradicts our choice of $\beta$.
Since $\Sk_{M_0}(A_{i+1})\subs \Sk_{M_{i+1}}(A_{i+1})$, we have
established the proposition.
\end{proof}

Note that the preceding proposition finishes the proof of the
claim
--- if play of the game continues for all $(2^\chi)^+$ steps, then
$\langle
\Sk_{M_0}(A_i)\cap\beth_{\alpha^*}(\chi):i<(2^\chi)^+\rangle$ is a
strictly decreasing family of subsets of $\Sk_{M_0}(A_{i^*})$,
contradicting (\ref{assumption}).
\end{proof}

\section{Club--guessing technology}

In this section, we prove that if $\lambda=\mu^+$, where $\mu$ is
singular, then under certain circumstances we can find a
complicated ``library'' of colorings of smaller cardinals.  In the
next section, we will use this library of colorings to get a
complicated coloring of $\lambda$.

The basics of club--guessing are explained in \cite{Sh:365}, but
we will take a few minutes to recall some of the definitions.

Let us recall that if $S$ is a stationary subset of $\lambda$,
then an $S$--club system is a sequence $\bar{C}=\langle
C_\delta:\delta\in S\rangle$ such that for (limit) $\delta\in S$,
$C_\delta$ is closed unbounded in $\delta$.

In this section, we will be concerned with the case where
$\lambda$ is the successor of a singular cardinal, i.e.,
$\lambda=\mu^+$ where $\cf(\mu)<\mu$. In this context, if
$\bar{C}$ is an $S$--club system, then for $\delta\in S$ we define
an ideal $J_\delta^{b[\mu]}$ on $C_\delta$ by $A\in
J_{\delta}^{b[\mu]}$ if and only if $A\subs C_\delta$, and for
some $\theta<\mu$ and $\gamma<\delta$,

\begin{equation*}
\beta\in A\cap\nacc(C_\delta)\Rightarrow\left[\beta<\gamma\text{
or }\cf(\beta)<\theta\right].
\end{equation*}
Note that it is a bit easier to understand the definition of
$J_\delta^{b[\mu]}$ by looking at the contrapositive --- a subset
$A$ of $C_\delta$ is ``large'', i.e., not in $J_\delta^{b[\mu]}$,
if and only if $A\cap\nacc(C_\delta)$ is cofinal in $\delta$, and
the cofinalities of members of any end segment of
$A\cap\nacc(C_\delta)$ are unbounded below $\mu$.

\begin{claim}
\label{1.5} Let $\lambda=\mu^+$, where $\mu$ is a singular
cardinal of cofinality $\kappa<\mu$. Let $S\subs\lambda$ be
stationary, and assume that $\sup\{\cf(\delta):\delta\in
S\}=\mu^*<\mu$. Let $\bar{C}$ be an $S$--club system, and for each
$\delta\in S$, let $J_\delta$ be the ideal $J_\delta^{b[\mu]}$.
Let $\langle\kappa_i:i<\kappa\rangle$ be a non--decreasing
sequence of cardinals such that
\begin{equation}
\kappa^*=\sum_{i<\kappa}\kappa_i\leq\mu,
\end{equation}
and let $\gamma^*<\mu$.

Assume we are given a $\lambda$--club system $\bar{e}$ and a
sequence of ideals $\bar{I}=\langle
I_\alpha:\alpha<\lambda\rangle$ such that
\begin{enumerate}
\item $I_\alpha$ is an ideal on $e_\alpha$ extending $J^{\bd}_{e_\alpha}$
\sk
\item \label{assumption1} if $\delta\in S$, then for each $i<\kappa$,
\begin{equation*}
\{\alpha\in\nacc(C_\delta):\text{ Player I wins
}\gm^\omega_{I_\alpha}[e_\alpha,\kappa_i,\gamma^*]\}=\nacc(C_\delta)\mod
J_\delta
\end{equation*}
\sk
\item \label{gengar} for any club $E\subs\lambda$, for stationarily many $\delta\in S$,
\begin{equation*}
\{\alpha\in\nacc(C_\delta):B_0[E,e_\alpha]\notin I_\alpha\}\notin
J_\delta,
\end{equation*}
where
\begin{equation*}
B_0[E, e_\alpha^*]=\{\beta\in\nacc(e_\alpha):E\text{ meets the
interval }(\sup(\beta\cap e_\alpha),\beta)\}.
\end{equation*}
\end{enumerate}
Then there is a function $h:\lambda\rightarrow(\kappa+1)$ and a
sequence
$$\bar{F}=\langle F_\delta:\delta<\lambda,\:\delta\text{ a limit }\rangle$$
such that
\begin{list}{$\circledast_1$}{\setlength{\leftmargin}{.5in}\setlength{\rightmargin}{.5in}}
\item
$F_\delta:[e_\delta]^{<\omega}\longrightarrow\kappa_{h(\delta)}\text{
(where $\kappa^*:=\kappa_\kappa$ )}$
\end{list}
and
\begin{list}{$\circledast_2$}{\setlength{\leftmargin}{.5in}\setlength{\rightmargin}{.5in}}
\item for every club $E\subs\lambda$, for each $i<\kappa$ there are stationarily many $\delta\in S$ such
 that the set of $\beta\in\nacc(C_\delta)$ satisfying the
 following
\begin{itemize}
\item $h(\beta)\geq i$
\sk
\item $B_0[E, e_\beta]\notin I_\beta$
\sk
\item for all $\gamma<\beta$, $\kappa_{h(\beta)}\subs\ran(F_\beta\restr\bigl{[}B_0[E, e_\beta]\setminus\gamma\bigr{]}^{<\omega})$
\end{itemize}
is {\em not} in $J_\delta$.
\end{list}
\end{claim}

Now admittedly the previous claim is quite a lot to digest, so we
will take a little time to illuminate the basic situation we have
in mind.

\begin{claim}
\label{helpclaim} The assumptions of Claim \ref{1.5} are satisfied
if
\begin{enumerate}
\item $\lambda=\mu^+$ where $\kappa=\cf(\mu)<\mu$
\sk
\item $S\subs \{\delta<\lambda:\cf(\delta)=\kappa\}$
\sk
\item $\delta\in S\rightarrow |\delta|=\mu$ (i.e., $S\subs \lambda\setminus \mu$)
\sk
\item $\bar{C}$ is an $S$--club system
\sk
\item $\bar{J}=\langle J_\delta:\delta\in S\rangle$ where $J_\delta=J^{b[\mu]}_{C_\delta}$
\sk
\item $\id_p(\bar{C},\bar{J})$ is a proper ideal
\sk
\item $\langle \kappa_i:i<\kappa\rangle$ is a non--decreasing sequence of cardinals with supremum $\kappa^*\leq\mu$
\sk
\item $\gamma^*<\mu$, and for each $i<\kappa$, Player I wins the game $\gm^\omega[\theta,\kappa_i,\gamma^*]$ for all large enough regular $\theta<\mu$
\sk
\item  $\bar{e}$ is a $\lambda$--club system such that $|e_\beta|<\mu$
\sk
\item for $\alpha<\lambda$, $I_\alpha=J^{\text{bd}}_{e_\alpha}$
\sk
\end{enumerate}
\end{claim}
\begin{proof}[Proof of Claim \ref{helpclaim}]
We need only check items (2) and (3) in the statement of Claim
\ref{1.5} --- everything else is trivially satisfied. Concerning
(2),  given $\delta\in S$ and $i<\kappa$, we need to show
\begin{equation*}
\{\alpha\in\nacc(C_\delta):\text{ Player I wins
}\gm^\omega[e_\alpha,\kappa_i,\gamma^*]\}=\nacc(C_\delta)\mod
J_\delta.
\end{equation*}
Let $A$ consist of those $\alpha\in\nacc(C_\delta)$ for which
Player I does not win the game
$\gm^{\omega}[e_\alpha,\kappa_i,\gamma^*]$. By our assumptions,
there is a $\theta<\mu$ such that $|e_\alpha|<\theta$ for all
$\alpha\in A$, and therefore $A$ is in the ideal
$J^{b[\mu]}_{C_\delta}=J_\delta$ and we have what we need.

Concerning (3), given $E\subs\lambda$ club, we must find
stationarily many $\delta\in S$ such that
\begin{equation*}
\{\alpha\in\nacc(C_\delta):B_0[E,e_\alpha]\notin I_\alpha\}\notin
J_\delta.
\end{equation*}
Let $E'=\{\xi\in E: \otp(E\cap\xi)=\xi\text{and $\mu$ divides
$\xi$}\}$. Clearly $E'$ is a closed unbounded subset of $E$, and
since $\id_p(\bar{C},\bar{J})$ is a proper ideal, the set
\begin{equation*}
S^*:=\{\delta\in S\cap E' : E'\cap \nacc(C_\delta)\notin
J_\delta\}
\end{equation*}
is stationary.

Fix $\delta\in S^*$, and suppose we are given $\theta<\mu$ and
$\xi<\delta$. Since $E'\cap\nacc(C_\delta)\notin J_\delta$, we can
find $\alpha\in E'\cap\nacc(C_\delta)$ such that
$\alpha>\max\{\xi,\mu\}$ and $\cf(\alpha)>\theta$.
 Since the order--type of $E\cap\alpha$ is
$\alpha\geq\mu>|e_\alpha|$, we know that $B_0[E, e_\alpha]$ is
unbounded in $e_\alpha$ hence a member of $I_\alpha$.  This shows
that the set of such $\alpha$ is in $J_\delta^+$, as required.
\end{proof}

Now we return to the proof of Claim \ref{1.5}.

\begin{proof}[Proof of Claim \ref{1.5}]
Let $\sigma=\cf(\sigma)$ be a regular cardinal $<\mu$ that is
greater than $\mu^*$ and $\gamma^*$.  For each limit
$\beta<\lambda$, if there is an $i\leq\kappa$ such that Player I
wins the version of
$\gm^\omega_{I_\beta}[e_\beta,\kappa_i,\sigma^+]$ where we allow
Player II to pass, then we let $h(\beta)$ be the maximal such $i$
--- note that $i$ exists by (\ref{usedgame})  of Claim \ref{1.4} ---
and let $\str_\beta$ be a strategy that witnesses this.

Note that since $\gamma^*<\sigma^+$ and
$J_\delta=J^{b[\mu]}_\delta$ for $\delta\in S$,  we have that for
$\delta\in S$ and $i<\kappa$ that
\begin{equation*}
\{\beta\in\nacc(C_\delta):\str_\beta\text{ is defined and }i\leq
h(\beta)\}=\nacc(C_\delta)\mod J_\delta.
\end{equation*}

We will make $\sigma^+$ attempts to build $\bar{F}$ witnessing the
conclusion. In stage $\zeta<\sigma^+$, we assume that our prior
work has furnished us with a decreasing sequence $\langle
E_\xi:\xi<\zeta\rangle$ of clubs in $\lambda$, and, for each
$\beta<\lambda$ where $\str_\beta$ is defined, an initial segment
$\langle F^\xi_\beta, A^\xi_\beta:\xi<\zeta\rangle$ of a play of
$\gm^\omega_{I_\beta}[e_\beta, \kappa^*_{h(\beta)},\sigma^+]$ in
which Player I uses $\str_\beta$. (Note that our convention is
that if Player II chooses to pass at a stage, we let $A^\xi_\beta$
be undefined.)

For each such $\beta$, let
$F^\zeta_\beta:[e_\beta]^{<\omega}\rightarrow\kappa_{h(\beta)}$ be
given by $\str_\beta$, and for those $\beta$ for which
$\str_\beta$ is undefined, we let $F^\beta_\zeta$ be any such
function. Now if $\langle
F^\zeta_\beta:\beta<\lambda\rangle:=\bar{F}^\zeta$ is as required
then we are done. Otherwise, there is a club $E'\subs\lambda$ and
$i_\zeta<\kappa$ exemplifying the failure of $\bar{F}^\zeta$, and
without loss of generality,
\begin{equation}
(\forall\delta\in S)\bigl{[}B_{i_\zeta}[E'_\zeta, C_\delta,
\bar{I},\bar{e},\bar{F}^\zeta]\bigr{]}\in J_\delta.
\end{equation}

Now let $E_\zeta=\acc(E'_\zeta\cap\bigcap_{\xi<\zeta}E_\xi)$.  For
each $\beta$ where $\str_\beta$ is defined, we let Player II
respond to $F^\zeta_\beta$ by playing the set
$B_0[E_\zeta,e_\beta]$ if it is a legal move, otherwise we let him
pass. We then proceed to stage $\zeta+1$.

Assuming that this construction continues for all $\sigma^+$
stages, we will arrive at a contradiction. Let
$E=\bigcap_{\zeta<\sigma^+} E_\zeta$. By assumption (\ref{gengar})
there is a $\delta(*)\in S$ for which
\begin{equation*}
A_1:=\{\beta\in\nacc(C_{\delta(*)}):B_0[E, e_\beta]\notin
I_\beta\}\notin J_{\delta(*)}.
\end{equation*}
By assumption (\ref{assumption1}), we have
\begin{equation*}
A_2:=\{\beta\in A_1: \str_\beta\text{ is defined }\}\notin
J_{\delta(*)}.
\end{equation*}

For $\beta\in A_2$, look at the play $\langle F_\beta^\zeta,
A^\zeta_\beta: \zeta<\sigma^+\rangle$. Since Player I wins, there
is a $\zeta_\beta<\sigma^+$ such that Player II passed at stage
$\zeta$ for all $\zeta\geq\zeta_\beta$.  Since $\sigma>\mu^*$ and
$J_{\delta(*)}$ is $\mu^*$--based, for some $\zeta^*<\sigma^+$,
\begin{equation*}
A_3=\{\beta\in A_1:\str_\beta\text{ is defined and
}\zeta_\beta\leq\zeta^*\}\notin J_{\delta(*)}.
\end{equation*}

Now $E_{\zeta^*}$ was defined so that for some $i_{\zeta^*}$, for
all $\delta\in S$,
\begin{equation}
\label{crud} B_{i_{\zeta^*}}[E_{\zeta^*}, C_{\delta},
\bar{I},\bar{e},\bar{F}^{\zeta^*}]\in J_\delta,
\end{equation}
but (again by assumption (\ref{assumption1}))
\begin{equation*}
A_4=\{\beta\in A_1:\str_{\beta}\text{ is defined, }
\zeta_\beta\leq\zeta^*,\text{ and }i_{\zeta^*}\leq
h(\beta)\}\notin J_{\delta(*)}.
\end{equation*}
 For $\beta\in A_4$, we know that at stage $\zeta^*$ of our play of
 $\gm^\omega_{I_\beta}[e_\beta, \kappa_{h(\beta)},\sigma^+]$ the set $B_0[E_{\zeta^*}, e_\beta]$
 was not a legal move.  Since our sequence of clubs is  decreasing, we know that
  $B_0[E_{\zeta^*}, e_\beta]$ is a subset of $B_0[E_\xi, e_\beta]$ for all $\xi<\zeta^*]$,
  so we have $$B_0[E_{\zeta^*}, e_\beta]\subs\bigcap_{\xi<\zeta^*} A^\xi_\beta.$$

Since $\beta\in A_1$, we know that $B_0[E_{\zeta^*},
e_\beta]\notin I_\beta$. Thus the reason for $B_0[E_{\zeta^*},
e_\beta]$ being an illegal move must be that for all
$\gamma<\beta$,
\begin{equation*}
\kappa^*_{h(\beta)}\subs\ran(F^{\zeta^*}_\beta\restr
[B_0[E_{\zeta^*}, e_\beta]\setminus\gamma]^{<\omega}).
\end{equation*}
All of these facts combine to tells us that $\beta\in
B_{i_{\zeta^*}}[E_{\zeta^*}, C_\delta, \bar{I},\bar{e},
\bar{F}^{\zeta^*}]$, and thus

\begin{equation*}
A_4\subs B_{i_{\zeta^*}}[E_{\zeta^*}, C_\delta, \bar{I},
\bar{e}^*,\bar{F}^{\zeta^*}]\notin J_{\delta(*)},
\end{equation*}
contradicting (\ref{crud}).
\end{proof}

The proofs in this section (and the next) can be considerably
simplified if we are willing to restrict ourselves to the case
 $\kappa^*<\mu$, as we can dispense with the sequence $\langle\kappa_i:i<\kappa\rangle$.

\section{Building the Coloring}

We now come to the main point of this paper; we dedicate this
section and the next to proving the following theorem.

\begin{thm}
\label{ash} Assume $\lambda=\mu^+$, where $\mu$ is a singular
cardinal of uncountable cofinality, say
$\aleph_0<\kappa=\cf(\mu)<\mu$. Assume $\langle
\kappa_i:i<\kappa\rangle$ is non--decreasing with supremum
$\kappa^*\leq\mu$, and there is a $\gamma^*<\mu$ such that
 for each $i$, for every large enough regular $\theta<\mu$,
Player I has a winning strategy in the game
$\gm^\omega[\theta,\kappa_i,\gamma^*]$. Then
$\pr_1(\lambda,\lambda,\kappa^*,\kappa)$ holds.
\end{thm}
Let $\langle S_i:i<\kappa\rangle$ be a sequence of pairwise
disjoint stationary subsets of
$\{\delta<\lambda:\cf(\delta)=\kappa\}$.  For $i<\kappa$, let
$\bar{C}^i$ be an $S_i$--club system such that 
\sk
\begin{itemize}
\item  $\lambda\notin\id_p(\bar{C}^i,\bar{J}^i)$, where
$\bar{J}^i=\langle J^{b[\mu]}_{C^i_\delta}:\delta\in S_i\rangle$
\sk
\item for $\delta\in S_i$, $\otp(C^i_\delta)=\cf(\delta)=\kappa=\cf(\mu)$
\end{itemize}
Such ladder systems can be found by Claim 2.6 (and Remark 2.6A
(6)) of \cite{Sh:g} --- for the second statement to hold, we need
that $\mu$ has uncountable cofinality.

\begin{claim}
There is a $\lambda$--club system $\bar{e}$ such that
$|e_\beta|\leq\cf(\beta)+\cf(\mu)$, and $\bar{e}$ ``swallows''
each $\bar{C}^i$, i.e., if $\delta\in
S_i\cap(e_\beta\cup\{\beta\})$, then $C^i_\delta\subs e_\beta$.
\end{claim}
\begin{proof}
Let $S=\cup_{i<\kappa}S_i$, and let $\beta<\lambda$ be a limit
ordinal.  Let $e^0_\beta$ be a closed cofinal subset of $\beta$ of
order--type $\cf(\beta)$. We will construct the required ladder
$e_\beta$ in $\omega$--stages, with $e_\beta^n$ denoting the
result of the first $n$ stages of our procedure. The construction
is straightforward, but it is worthwhile to note that we need to
use the fact that each member of $S$ has uncountable cofinality.

Given $e^n_\beta$, let us define
\begin{equation}
B_n=S\cap (e^n_\beta\cup\{\beta\}).
\end{equation}

Now we let $e^{n+1}_\beta$ be the closure in $\beta$ of
\begin{equation}
e^{n}_\beta\cup\bigcup\{C_\delta:\delta\in B_n\}.
\end{equation}

Note that $|e^{n+1}|\leq \cf(\mu)+\cf(\beta)$ as
$|C_\delta|=\cf(\mu)=\kappa$ for each $\delta\in S$.  Finally, we
let $e_\beta$ be the closure of $\cup_{n<\omega}e^n_\beta$ in
$\beta$.

Clearly $|e_\beta|\leq\cf(\mu)+\cf(\beta)$. Also, since each
element of $S$ has uncountable cofinality, if $\delta\in S\cap
e_\beta$, then there is an $n$ such that $\delta\in e^n_\beta$,
and therefore
\begin{equation}
C_\delta\subs e^{n+1}_\beta\subs e_\beta,
\end{equation}
as required.
\end{proof}
\sk

For each $i<\kappa$, there are $h_i$ and $\bar{F}^i=\langle
F^i_\delta:\delta<\lambda,\; \delta\text{ limit }\rangle$ as in
the conclusion of Claim \ref{1.5} applied to $\bar{C}^i$ and
$\bar{e}$; note that we satisfy the assumptions of Claim
\ref{1.5} by way of Claim \ref{helpclaim}.

Let $\langle\lambda_i:i<\kappa\rangle$ be a strictly increasing
sequence of regular cardinals $>\kappa$ and cofinal in $\mu$ such
that
\begin{equation}
\lambda=\tcf\bigl{(}\prod_{i<\kappa}\lambda_i/J^{\bd}_\kappa\bigr{)},
\end{equation}
and let $\langle f_\alpha:\alpha<\lambda\rangle$ exemplify this.
Finally, let $h_0^*:\kappa\rightarrow\omega$ and
$h_1^*:\kappa\rightarrow\kappa$ be such that
\begin{equation}
(\forall n)(\forall i<\kappa)(\exists^\kappa
j<\kappa)[h^*_0(j)=n\text{ and }h_1^*(j)=i].
\end{equation}

Before we can define our coloring, we must recall some of the
terminology of \cite{Sh:g}.

\begin{defn}
Let $0<\alpha<\beta<\lambda$, and define
\begin{equation*}
\gamma(\alpha,\beta)=\min\{\gamma\in e_\beta:\gamma\geq\alpha\}.
\end{equation*}
We also define (by induction on $\ell$)
\begin{equation*}
\gamma_0(\alpha,\beta)=\beta,
\end{equation*}
\begin{equation*}
\gamma_{\ell+1}(\alpha,\beta)=\gamma(\alpha,\gamma_\ell(\alpha,\beta))\text{
(if defined).}
\end{equation*}
We let $k(\alpha,\beta)$ be the first $\ell$ for which
$\gamma_\ell(\alpha,\beta)=\alpha$.  The sequence $\langle \gamma_i(\alpha,\beta):i\leq k(\alpha,\beta)\rangle$ will be referred to as the {\em walk from $\beta$ to $\alpha$ along the ladder system $\bar{e}$}.
\end{defn}

We now define the coloring $c$ that will witness
$\pr_1(\lambda,\lambda,\kappa^*,\kappa)$. Recall that $c$ must be
a symmetric two--place function from $\lambda$ to $\kappa^*$.

Given $\alpha<\beta$, we let $i=i(\alpha,\beta)$ be the maximal
$j<\kappa$ such that $f_\beta(j)<f_\alpha(j)$ (if such an $j$
exists).  Next, we walk from $\beta$ down to $\alpha$ along
$\bar{e}$ until we reach an ordinal $\nu(\alpha,\beta)$ such that
\begin{equation*}
f_\alpha(i)<f_{\nu(\alpha,\beta)}(i),
\end{equation*}
(again, if such an ordinal exists.)  After this, we walk along $\bar{e}$ from
$\alpha$ toward the ordinal $\max(\alpha\cap
e_{\nu(\alpha,\beta)})$ until we reach an ordinal
$\eta(\alpha,\beta)$ for which
\begin{equation*}
f_{\nu(\alpha,\beta)}(i)<f_{\eta(\alpha,\beta)}(i).
\end{equation*}

The idea now is to look at how the ladders $e_{\nu(\alpha,\beta)}$
and $e_{\eta(\alpha,\beta)}$ intertwine.  Let us make a temporary
definition by calling an ordinal $\xi\in e_{\nu(\alpha,\beta)}$
{\em relevant} if $e_{\eta(\alpha,\beta)}$ meets the interval
$(\sup(\xi\cap e_{\nu(\alpha,\beta)}), \xi)$.

If it makes sense, we let $w(\alpha,\beta)\subs
e_{\nu(\alpha,\beta)}$ be the last $h_0^*(i(\alpha,\beta))$
relevant ordinals in $e_{\nu(\alpha,\beta)}$ (so we need that the
relevant ordinals have order--type $\gamma+h_0^*(i(\alpha,\beta))$
for some $\gamma$).

Finally, we define our coloring by
\begin{equation}
c(\alpha,\beta)=F_{\nu(\alpha,\beta)}^{h_1^*(i(\alpha,\beta))}(w(\alpha,\beta)).
\end{equation}

If the attempt to define $c(\alpha,\beta)$ breaks down at some
point for some specific $\alpha<\beta$, then we set
$c(\alpha,\beta)=0$.

We now prove that this coloring works, so suppose $\langle
t_\alpha:\alpha<\lambda\rangle$ are pairwise disjoint subsets of
$\lambda$ such that $|t_\alpha|=\theta_1<\kappa$ and
$j^*<\kappa^*$, and without loss of generality $\alpha<\min
t_\alpha$ and $\theta_1\geq\omega$. We need to find $\delta_0$ and
$\delta_1$ such that
\begin{equation}
 \alpha\in t_{\delta_0}\text{ and }\beta\in t_{\delta_1}
 \Rightarrow \alpha<\beta\text{ and }c(\alpha,\beta)=j^*.
\end{equation}

Let $j_1$ be the least $j$ such that $j^*<\kappa_j$, and let $S$,
$\bar{C}$, and $\bar{F}$ denote $S_{j_1}$, $\bar{C}^{j_1}$, and
$\bar{F}^{j_1}$ respectively.

Given $\delta<\lambda$, we define the {\em envelope of $t_\delta$}
(denoted $\env(t_\delta)$) by the formula
\begin{equation}
\env(t_\delta)=\bigcup_{\zeta\in t_\delta}\{\gamma_\ell(\delta,\zeta):\ell\leq k(\delta,\zeta)\}.
\end{equation}
The envelope of $t_\delta$ is the set of all ordinals obtained by
walking down to $\delta$ from some $\zeta\in t_\delta$ using the
ladder system $\bar{e}$. This makes sense as we have arranged that
$\delta<\min t_\delta$. Note also that $|\env(t_\delta)|\leq
|t_\delta|=\theta_1$.

Next we define functions $g^{\min}_\delta$ and $g^{\max}_\delta$
in $\prod_{i<\kappa}\lambda_i$ by
\begin{equation}
g^{\min}_\delta(i)=\min\{f_\gamma(i):\gamma\in\env(t_\delta)\},
\end{equation}
and
\begin{equation}
g^{\max}_\delta(i)=\sup\{f_\gamma(i)+1:\gamma\in\env(t_\delta)\}.
\end{equation}
Note that $g^{\max}_\delta$ is well--defined as we assume that
$\kappa<\min\{\lambda_i:i<\kappa\}$.

The following claim is quite easy, and the proof is left to the
reader.

\begin{claim}\hfill
\begin{enumerate}
\item $f_\delta=_{\jbk} g^{\min}_\delta$
\sk
\item $g^{\min}_\delta(i)\leq g^{\max}_\delta(i)$ for all $i<\kappa$
\sk
\item There is a $\delta'>\delta$ such that $g^{\max}_\delta\leq_{\jbk} g^{\min}_{\delta'}$.
\sk
\end{enumerate}
\end{claim}

Now let $\chi^*=(2^\lambda)^+$, and let $\langle
M_\alpha:\alpha<\lambda\rangle$ be a sequence of elementary
submodels of $\langle H(\chi^*),\in,<^*_{\chi^*}\rangle$  that is
increasing and continuous in $\alpha$ and such that each
$M_\alpha\cap\lambda$ is an ordinal, $\langle
M_\beta:\beta\leq\alpha\rangle\in M_{\alpha+1}$, and $\langle
f_\alpha:\alpha<\lambda\rangle$, $g$, $c$, $\bar{e}$, $S$,
$\bar{C}$, $\langle t_\alpha:\alpha<\lambda\rangle$ all belong to
$M_0$. Note that $\mu+1\subs M_0$.

The set $E=\{\alpha<\lambda: M_\alpha\cap\lambda=\alpha\}$ is
closed unbounded in $\lambda$, and furthermore,
\begin{equation}
\alpha<\delta\in E\Rightarrow \sup t_\alpha<\delta.
\end{equation}

By the choice of $\bar{C}$ and $\bar{F}$, for some $\delta\in
S\cap E$ we have the set
\begin{equation}
A=\{\beta\in\nacc(C_\delta):(\forall\gamma<\beta)\ran(F_\beta\restr\bigl{[}B_0[E,
e_\beta]\setminus\gamma\bigr{]}^{<\omega})\sups\kappa_{j_1}\}
\end{equation}
is not in $J^{b[\mu]}_{C_\delta}$.

Note that $A\subs\acc(E)$, as $B_0[E,e_\beta]$ is unbounded in
$\beta$ for $\beta\in A$. For $\beta\in t_\delta$, if
$\ell<k(\delta,\beta)$ then
$e_{\gamma_\ell(\delta,\beta)}\cap\delta$ is bounded in $\delta$,
and since it is closed it has a well--defined maximum. Since
$|t_\delta|<\kappa=\cf(\delta)$, this means the ordinal
\begin{equation*}
\gamma^{\otimes}:=\sup\{\max[e_{\gamma_\ell(\delta,\beta)}\cap\delta]:\beta\in
t_\delta\text{ and }\ell<k(\delta,\beta)\}
\end{equation*}
is strictly less than $\delta$.

For $\beta\in t_\delta$, let us define
\begin{equation}
A_\beta:=\{\beta'\in A: (\exists\ell\leq
k(\beta,\delta))[\cf(\beta')\leq
|e_{\gamma_\ell(\delta,\beta)}|]\}.
\end{equation}
Since the cardinality of each ladder in $\bar{e}$ is less than
$\mu$, each set $A_\beta$ is an element of
$J^{b[\mu]}_{C_\delta}$.  The ideal $J^{b[\mu]}_{C_\delta}$ is
$\kappa$--complete, so the fact that $|t_\delta|<\kappa$ and
$k(\beta,\delta)$ is finite for each $\beta\in t_\delta$ together
imply that
\begin{equation}
\bigcup_{\beta\in t_\delta} A_\beta\in J^{b[\mu]}_{C_\delta}.
\end{equation}
By the definition of $A$ and our choice of $\delta$, this means it
is possible to choose $\beta^*\in A\setminus(\gamma^{\otimes}+1)$
that is not in any $A_\beta$, i.e.,
\begin{equation}
\beta\in t_\delta\text{ and }\ell< k(\delta,\beta)\Longrightarrow
\cf(\beta^*)>|e_{\gamma_\ell(\delta,\beta)}|.
\end{equation}

\begin{claim}\hfill
\label{pikachu}
\begin{enumerate}
\item If $\epsilon\in t_\delta$, and $\ell=k(\delta,\epsilon)-1$, then
 $\beta^*\in\nacc(e_{\gamma_\ell(\delta,\epsilon)})$.
\sk
\item If $\epsilon\in t_\delta$ and $\gamma^{\otimes}<\gamma'\leq\beta^*$, then
\sk
\begin{itemize}
\item $\gamma_\ell(\delta,\epsilon)=\gamma_\ell(\gamma',\epsilon)$ for $\ell<k(\delta,\epsilon)$, and
\sk
\item $\gamma_{k(\delta,\epsilon)}(\gamma',\epsilon)=\beta^*$
\sk
\end{itemize}
\end{enumerate}
\end{claim}

\begin{proof}
For the first clause, note that $\delta$ is an element of
$e_{\gamma_\ell(\delta,\epsilon)}$ and hence by our choice of
$\bar{e}$, $C_\delta\subs e_{\gamma_{\ell}(\delta,\epsilon)} $.
Thus $\beta^*\in e_{\gamma_{\ell}(\delta,\epsilon)}$, and since
$\cf(\beta^*)>|e_{\gamma_{\ell}(\delta,\epsilon)}|$, we know that
$\beta^*$ cannot be an accumulation point of
$e_{\gamma_{\ell}(\delta,\epsilon)}$.

The first part of the second statement follows because of the
definition of $\gamma^{\otimes}$.  As far as the second part of
the second statement goes, it is best visualized as follows:

We walk down the ladder system $\bar{e}$ from $\epsilon$ to
$\gamma'$, we eventually hit a ladder that contains $\delta$ ---
this happens at stage $k(\delta,\epsilon)-1$. Since $C_\delta$ is
a subset of this ladder, the next step in our walk from $\epsilon$
to $\gamma'$ must be down to $\beta^*$ because
$\gamma^{\otimes}<\gamma'<\beta^*$.
 \end{proof}

We can visualize the preceding claim in the following manner:
$\beta^*$ is chosen so that for all sufficiently large
$\gamma'<\beta^*$,  all the walks from some element of $t_\delta$
to $\gamma'$ are funnelled through $\beta^*$
--- $\beta^*$ acts as a bottleneck. This will be key when want to prove that
our coloring works.

Since $\beta^*\in A$,  we can choose a finite increasing sequence
$\xi_0<\xi_1<\dots<\xi_n$ of ordinals in
$\acc(E)\cap\nacc(e_{\beta^*})\setminus(\gamma^\otimes+1)$ such
that $F^{j_1}_{\beta^*}(\{\xi_0,\dots,\xi_n\})=j^*$, the color we
are aiming for.

For each $\ell\leq n$, we can find $\zeta_\ell\in
E\setminus(\gamma^\otimes+1)$ such that
$$\sup(e_{\beta^*}\cap\xi_\ell)<\zeta_\ell<\xi_\ell.$$

Now we let $\phi(x_0,y_0,x_1,y_1,\dots, x_n, y_n, z_0, z_1)$ be
the formula (with parameters $\gamma^\otimes$, $\bar{f}$,
$\langle\lambda_i:i<\kappa\rangle$, $\bar{C}$, $\bar{e}$, $\langle
t_\alpha:\alpha<\lambda\rangle$, $h$, $h_0$, $j^*$) that describes
our current situation with $x_\ell$, $y_\ell$ standing for $\zeta_\ell$, $\xi_\ell$, and $z_0$, $z_1$ standing for $\beta^*$, $\delta$ , i.e., $\phi$ states
\begin{itemize}
\item $\gamma^\otimes<x_0<y_0<\dots<x_n<y_n<z_0<z_1$ are ordinals $<\lambda$
\sk
\item $z_1\in S$ and $z_0\in\nacc(C_{z_1})$
\sk
\item $\gamma^\otimes=\sup\{\max[e_{\gamma_\ell(z_1, \zeta)}\cap z_1]:\ell<k(z_1,\zeta)\text{ and }\zeta\in t_{z_1}\}$
\sk
\item $z_0\in\nacc(e_{\gamma_{k(z_1,\epsilon)}(z_1,\epsilon)})$ for all $\epsilon\in t_{z_1}$
\sk
\item $F^{j_1}_{z_0}(\{y_0,\dots,y_n\})=j^*$
\sk
\end{itemize}
Now clearly we have
\begin{equation}
\model\models\phi[\zeta_0,\xi_0,\dots,\zeta_n,\xi_n,\beta^*,\delta].
\end{equation}
%%%%%%%%%%%%%%
%
%   z_0 is \beta^*
%
%   z_1 is \delta
%
%
%
%%%%%%%%%%%%%
Recall that all the parameters needed in $\phi$ are in $M_0$, except possibly for $\gamma^\otimes$, so the model $M_{\gamma^\otimes+1}$ contains all the parameters we need.  Also, 
$\{\zeta_0,\xi_0,\dots,\zeta_n,\xi_n\}\in M_{\beta^*}$,
$\beta^*\in M_\delta\setminus M_{\beta^*}$, and since $\delta\in
\lambda\setminus M_\delta$, we have (recalling that $\exists^*z<\lambda$ means ``for unboundedly many $z<\lambda$)

\begin{equation}
M_\delta\models (\exists^*
z_1<\lambda)\phi(\zeta_0,\xi_0,\dots,\zeta_n,\xi_n,\beta^*, z_1).
\end{equation}
Therefore, this formula is true in $\model$ because of
elementarity. Similarly, we have
\begin{equation*}
\model\models(\exists^*z_0<\lambda)(\exists^*z_1<\lambda)\phi(\zeta_0,\xi_0,\dots,\zeta_n,\xi_n,z_0,z_1).
\end{equation*}
Now each of the intervals $[\gamma^\otimes+1,\zeta_0)$, $[\zeta_0,\xi_0)$,
$\dots$, contains a member of $E$,  so (by the definition of $E$)
similar considerations give us
\begin{equation*}
\model\models(\exists^*x_0<\lambda)\dots(\exists^*y_n<\lambda)(\exists^*z_0<\lambda)(\exists^*z_1<\lambda)\phi(x_0,y_0,\dots,z_0,z_1).
\end{equation*}

Now we can choose (in order)
\begin{equation}
\zeta^a_0<\zeta^b_0<\xi^a_0<\zeta^a_1<\xi^b_0<\zeta^b_1<\dots<\zeta^a_n
<\xi^b_{n-1}<\zeta^b_n<\xi^a_n
\end{equation}
such that
\begin{equation}
\label{venusaur}
(\exists^*z_0<\lambda)(\exists^*z_1<\lambda)[\phi(\zeta^a_0,\dots,\xi^a_{n-1},\zeta^a_n,\xi^a_n,z_0,z_1)],
\end{equation}
and
\begin{equation}
(\exists^*y_n<\lambda)(\exists^*z_0<\lambda)(\exists^*z_1<\lambda)[\phi(\zeta^b_0,\dots,\xi^b_{n-1},\zeta^b_n,y_n,z_0,z_1)],
\end{equation}

Our goal is to show that for all sufficiently large $i<\kappa$, it
is possible to choose objects $\beta^a$, $\delta^a$, $\xi^b_n$,
$\beta^b$, and $\delta^b$ such that \sk \sk

\begin{center}
\begin{tabular}{|l|}
\hline
\\
(1) $\zeta^b_{n}<\beta^a<\delta^a<\min(t_{\delta^a})\leq\max(t_{\delta^a})<\xi^b_n<\beta^b<\delta^b$\\
\\
(2) $\phi(\zeta^a_0,\dots,\xi^a_n,\beta^a,\delta^a)$\\
\\
(3) $\phi(\zeta^b_0,\dots,\xi^b_n,\beta^b,\delta^b)$\\
\\
(4) for all $\epsilon\in\env(t_{\delta^a})$, $g^{\min}_{\delta^a}\restr[i,\kappa)\leq f_\epsilon\restr[i,\kappa)\leq g^{\max}_{\delta^a}\restr [i,\kappa)$\\
\\
(5) for all $\epsilon\in\env(t_{\delta^b})$, $g^{\min}_{\delta^b}\restr[i,\kappa)\leq f_\epsilon\restr[i,\kappa)\leq g^{\max}_{\delta^b}\restr[i,\kappa)$\\
\\
(6) $g^{\max}_{\delta^b}(i)<g^{\min}_{\delta^a}(i)\leq g^{\max}_{\delta^a}(i)<f_{\beta^b}(i)< f_{\beta^a}(i)$\\
\\
(7)  $g^{\max}_{\delta^a}\restr[i+1,\kappa)<g^{\min}_{\delta^b}\restr[i+1,\kappa)$\\
\\
\hline
\end{tabular}
\end{center}
\sk
\begin{center}
\bf{Table 1}
\end{center}
\sk
\begin{claim}
\label{final} If for all sufficiently large $i<\kappa$ it is
possible to find objects satisfying the requirements of Table 1,
then we can find $\delta^a<\delta^b$ such that
$c(\epsilon^a,\epsilon^b)=j^*$ for all $\epsilon^a\in
t_{\delta^a}$ and $\epsilon^b\in t_{\delta^b}$.
\end{claim}
\begin{proof}
Let us choose $i^*<\kappa$ such that
\begin{itemize}
\item suitable objects (as above) can be found, and
\sk
\item  $h_1^*(i^*)=j_1$ and $h^*_0(i^*)=n$
\sk

\end{itemize}

Choose $\epsilon^a\in t_{\delta^a}$ and $\epsilon^b\in
t_{\delta^b}$; we verify that $c(\epsilon^a,\epsilon^b)=j^*$.

\begin{subclaim}
$i(\epsilon^a,\epsilon^b)=i^*$.
\end{subclaim}
\begin{proof}
Immediate by (4)-(7) in the table.
\end{proof}
\sk
\begin{subclaim}
$\nu(\epsilon^a,\epsilon^b)=\beta^b$.
\end{subclaim}
\begin{proof}
Note that $\gamma^\otimes<\epsilon^a<\beta^b$. Clause (3) of the table implies that the assumptions of Claim \ref{pikachu} hold. Thus by Claim
\ref{pikachu}, for $\ell<k(\delta^b,\epsilon^b)$ we have
\begin{equation*}
\gamma_\ell(\epsilon^a,\epsilon^b)=\gamma_\ell(\delta^b,\epsilon^b),
\end{equation*}
hence $\gamma_\ell(\epsilon^a,\epsilon^b)\in\env(t_{\delta^b})$
and (by (6) of the table and the definitions involved)
\begin{equation}
f_{\gamma_{\ell}(\epsilon^a,\epsilon^b)}(i^*)\leq
g^{\max}_{\delta^b}(i^*)<g^{\min}_{\delta^a}(i^*)\leq
f_{\epsilon^a}(i^*).
\end{equation}

For $\ell=k(\delta^b,\epsilon^b)$, Claim \ref{pikachu} tells us
\begin{equation*}
\gamma_\ell(\epsilon^a,\epsilon^b)=\beta^b,
\end{equation*}
and we have arranged that
\begin{equation}
f_{\epsilon^a}(i^*)\leq g^{\max}_{\delta^a}(i^*)<f_{\beta^b}(i^*).
\end{equation}
This establishes $\beta^b=\nu(\epsilon^a,\epsilon^b)$.
\end{proof}
\sk
\begin{subclaim}
$\eta(\epsilon^a,\epsilon^b)=\beta^a$.
\end{subclaim}
\begin{proof}
Let $\alpha=\max(e_{\beta^b}\cap\epsilon^a)$.  We have arranged that
\begin{equation*}
\zeta^b_n<\beta^a<\delta^a<\epsilon^a<\xi^b_n
\end{equation*}
and $\gamma^\otimes<\max(e_{\beta^b}\cap\delta^a)$, hence
$\gamma^\otimes<\alpha<\beta^a$.  For
$\ell<k(\delta^a,\epsilon^a)$, Claim \ref{pikachu} implies
\begin{equation*}
\gamma_\ell(\alpha,\epsilon^a)=\gamma_\ell(\delta^a,\epsilon^a)\in\env(t_{\delta^a}).
\end{equation*}
By our choice of $i^*$, we have
\begin{equation}
f_{\gamma_\ell(\alpha,\epsilon^a)}(i^*)\leq
g^{\max}_{\delta^a}(i^*)<f_{{\beta^b}}(i^*).
\end{equation}

For $\ell=k(\delta^a,\epsilon^a)$, Claim \ref{pikachu} implies
$\gamma_\ell(\alpha,\epsilon^a)=\beta^a$, and we have ensured
\begin{equation}
f_{\beta^b}(i^*)<f_{\beta^a}(i^*).
\end{equation}
Thus $\beta^a$ is the first ordinal $\eta$ in the walk from
$\epsilon^a$ to $\max(e_{\beta^b}\cap\epsilon^a)$ for which
$f_{\eta}(i^*)>f_{\beta^b}(i^*)$, and therefore
$\eta(\epsilon^a,\epsilon^b)=\beta^a$.
\end{proof}

\begin{subclaim}
$w(\epsilon^a,\epsilon^b)=\{\xi^b_0,\dots\xi^b_n\}$.
\end{subclaim}
\begin{proof}
Our previous subclaims imply that an ordinal $\xi\in e_{\beta^b}$
is relevant if and only if the ladder $e_{\beta^a}$ meets the
interval $(\sup(e_{\beta^b}\cap \xi),\xi)$.  Since
$h_0^*(i^*)=n+1$, we know that $w(\epsilon^a,\epsilon^b)$ consists
of the last $n+1$ relevant ordinals in $e_{\beta^b}$.

For $i\leq n$, clearly $\xi^b_i\in e_{\beta^b}$ and
$\sup(\xi^b_i\cap e_{\beta^b})\leq\zeta^b_n$.  We have made sure that
$e_{\beta^a}\cap (\zeta^b_i,\xi^b_i)\neq \emptyset$ (for example,
$\xi^a_{i}$ is an element in this intersection) and so each
$\xi^b_i$ is relevant.

Since $\beta^a<\xi^b_n$, it is clear that there are no relevant
ordinals larger than $\xi^b_n$.

Given $i<n$, if $\xi\in e_{\beta^b}\cap (\xi^b_i,\xi^b_{i+1})$,
then
\begin{equation*}
\xi^b_i\leq \sup(\xi\cap e_{\beta^b})\leq\xi\leq\zeta^b_{i+1}.
\end{equation*}
Since $\zeta^a_{i+1}<\xi^b_i<\zeta^b_{i+1}<\xi^{a}_{i+1}$, it
follows that
\begin{equation*}
[\sup(\xi\cap e_{\beta^b}),\xi)\subs [\zeta^a_{i+1},\xi^a_{i+1}),
\end{equation*}
and so $\xi$ is not relevant.  Thus $\{\xi^b_0,\dots,\xi^b_n\}$
are the last $n+1$ relevant elements of $e_{\beta^b}$, as was
required.
\end{proof}

To finish the proof of Claim \ref{final}, we note that as
$h_1^*(i^*)=j^*$, we have
\begin{equation}
c(\epsilon^a,\epsilon^b)=F^{j_1}_{\beta^b}(\{\xi^b_0,\dots,\xi^b_n\})=j^*.
\end{equation}
\end{proof}
\section{Finding the required ordinals}

The whole of this section will be occupied with showing that for
all sufficiently large $i<\kappa$, it is possible to find objects
satisfying the requirements of Table 1.

We begin with some notation intended to simplify the presentation.
\sk
\begin{itemize}
\item $\phi^a(z_0, z_1)$ abbreviates the formula $\phi(\zeta^a_0,\dots,\xi^a_n, z_0, z_1)$
\sk
\item $\phi^b(y_n, z_0, z_1)$ abbreviates the formula $\phi(\zeta_0^b,\zeta^n_b,y_n, z_0, z_1)$
\sk
\item For $i<\kappa$,  $\psi(i,z_1)$ abbreviates the formula
\begin{equation}
(\forall\epsilon\in
\env(t_{z_1}))[g^{\min}_{z_1}\restr[i,\kappa)\leq f_\epsilon\restr
[i,\kappa)\leq g^{\max}_{z_1}\restr[i,\kappa)]
\end{equation}
\end{itemize}

We have arranged things so that the sentence
\begin{multline}
(\exists^* z^a_0<\lambda)(\exists^* z_1^a<\lambda)(\exists^*
y_n^b<\lambda)\\(\exists^* z_0^b<\lambda)(\exists^*
z_1^b<\lambda)[\phi^a(z_0^a,z^a_1)\wedge\phi^b(y^b_n,z^b_0,z^b_1)]
\end{multline}
holds.

There are far too many alternations of quantifiers in the above
formula for most people to deal with comfortably; the best way to
view them is as a single quantifier that asserts the existence of
a tree of $5$--tuples with the property that every node of the
tree has $\lambda$ successors, and every branch through the tree
gives us five objects satisfying $\phi^a\wedge\phi^b$.

Let $\Phi(i,z^a_0,\dots,z^b_1)$ abbreviate the formula
\begin{multline*}
\phi^a(z^a_0,z^a_1)\wedge\phi^b(y^b_n, z^b_0, z^b_1)\wedge\psi(i,z^a_1)\wedge\psi(i,z^b_1)\\
\wedge\left(g^{\max}_{z^a_1}\restr[i+1,\kappa)<g^{\min}_{z^b_1}\restr[i+1,\kappa)\right).
\end{multline*}
By pruning the tree so that every branch through it is a strictly
increasing $5$--tuple, we get
\begin{multline}
\label{eqnyuck}
(\exists^* z^a_0<\lambda)(\exists^* z_1^a<\lambda)(\exists^* y_n^b<\lambda)\\
(\exists^* z_0^b<\lambda)(\exists^*
z_1^b<\lambda)(\forall^*i<\kappa)[\Phi(i,z^a_0,\dots,z^b_1)].
\end{multline}

We now make a rather {\em ad hoc} definition of another quantifier
in an attempt to make the arguments that follow a little bit
clearer. Given $i<\kappa$, let the quantifier
$\exists^{*,i}z^b_0<\lambda$ mean that not only are there
unboundedly many $z^b_0$'s below $\lambda$ satisfying whatever
property, but also that for each $\alpha<\lambda_i$, we can find
unboundedly many suitable $z^b_0$'s for which $f_{z^b_0}(i)$ is
greater than $\alpha$.

\begin{claim}
If we choose $\beta^a<\delta^a<\xi^b_n$ such that
\begin{equation}
\label{eqn1} (\exists^{*} z_0^b<\lambda)(\exists^*
z_1^b<\lambda)(\forall^*i<\kappa)[\Phi(i,\beta^a,\delta^a,\xi^b_n,z^b_0,z^b_1)],
\end{equation}
then
\begin{equation}
\label{eqn2} (\forall^*i<\kappa)(\exists^{*,i}
z_0^b<\lambda)(\exists^*
z_1^b<\lambda)[\Phi(i,\beta^a,\delta^a,\xi^b_n,z^b_0,z^b_1)].
\end{equation}
\end{claim}
\begin{proof}
Suppose that we have $\beta^a<\delta^a<\xi^b_n$ such that
(\ref{eqn1}) holds but (\ref{eqn2}) fails.  Then there is an
unbounded $I\subs\kappa$  such that for each $i\in I$,
\begin{equation}
\label{eqn3} \neg(\exists^{*,i} z_0^b<\lambda)(\exists^*
z_1^b<\lambda)[\Phi(i,\beta^a,\delta^a,\xi^b_n,z^b_0,z^b_1)].
\end{equation}
In (\ref{eqn1}), we can move the quantifier
``$\forall^*i<\kappa''$ past the quantifiers to its left, i.e.,
\begin{equation}
(\forall^*i<\kappa)(\exists^{*} z_0^b<\lambda)(\exists^*
z_1^b<\lambda)[\Phi(i,\beta^a,\delta^a,\xi^b_n,z^b_0,z^b_1)],
\end{equation}
so without loss of generality, for all $i\in I$,
\begin{equation}
(\exists^{*} z_0^b<\lambda)(\exists^*
z_1^b<\lambda)[\Phi(i,\beta^a,\delta^a,\xi^b_n,z^b_0,z^b_1)].
\end{equation}
Since (\ref{eqn3}) holds for all $i\in I$, it must be the case
that for each $i\in I$, there is a value $g(i)<\lambda_i$ such
that for all sufficiently large $\beta<\lambda$, if
\begin{equation}
\label{eqn4}
(\exists^*z^b_1<\lambda)[\Phi(i,\beta^a,\delta^a,\xi^b_n,\beta,z^b_1)],
\end{equation}
then
\begin{equation}
\label{eqn5} f_\beta(i)\leq g(i).
\end{equation}
Since $\{f_\alpha:\alpha<\lambda\}$ witnesses that the true
cofinality of $\prod_{i<\kappa}\lambda_i$ is $\lambda$, we know
\begin{equation}
(\forall^*x<\lambda)(\forall^*i\in I)[g(i)<f_x(i)].
\end{equation}
When we combine this with (\ref{eqn1}), we see that it is possible
to choose $\beta^b<\lambda$ such that
\begin{equation}
(\forall^*i\in I)[g(i)<f_{\beta^b}(i)],
\end{equation}
and
\begin{equation}
(\exists^*z^b_1<\lambda)(\forall^*j<\kappa)[\Phi(j,\beta^a,\delta^a,\xi^b_n,\beta^b,z^b_1)].
\end{equation}
(Note that we have quietly used the fact that
$|I|<\lambda=\cf(\lambda)$ to get a $\beta^b$ that is ``large
enough'' so that (\ref{eqn4}) implies (\ref{eqn5}) for  all $i\in
I$ for this particular $\beta^b$.) This last equation implies
\begin{equation*}
(\forall^*j<\kappa)(\exists^*z^b_1<\lambda)[\Phi(j,\beta^a,\delta^a,\xi^b_n,\beta^b,z^b_1)],
\end{equation*}
so it is possible to choose $i\in I$ large enough so that
\begin{equation*}
g(i)<f_{\beta^b}(i)
\end{equation*}
and
\begin{equation*}
(\exists^*z^b_1<\lambda)[\Phi(i,\beta^a,\delta^a,\xi^b_n,\beta^b,z^b_1)].
\end{equation*}
This is a contradiction, as (\ref{eqn4}) holds for our choice of
$i$ and $\beta=\beta^b$, yet (\ref{eqn5}) fails.
\end{proof}

Notice that an immediate corollary of the preceding claim is
\begin{multline}
(\exists^*z^a_0<\lambda)(\exists^*z^a_1<\lambda)(\exists^*y^b_n<\lambda)(\forall^*i<\kappa)\\
(\exists^{*,i}z^b_0<\lambda)(\exists^*
z_1^b<\lambda)[\Phi(i,\beta^a,\delta^a,\xi^b_n,z^b_0,z^b_1)].
\end{multline}

\begin{claim}
If $\beta^a<\lambda$ is chosen so that
\begin{multline}
\label{eqn6}
(\exists^*z^a_1<\lambda)(\exists^*y^b_n<\lambda)(\forall^*i<\kappa)\\
(\exists^{*,i}z^b_0<\lambda)(\exists^*
z_1^b<\lambda)[\Phi(i,\beta^a,z^a_1,y^b_n,z^b_0,z^b_1)],
\end{multline}
then
\begin{equation*}
(\forall^*i<\kappa)(\exists v<\lambda_i)(\exists^*z^a_1<\lambda)[\psi'\wedge\psi'']
\end{equation*}
where
\begin{equation*}
\psi':=g^{\max}_{z^a_1}(i)<v,
\end{equation*}
and
\begin{equation*}
\psi'':=(\exists^*y^b_n<\lambda)(\exists^*z^b_0<\lambda)\left[v<f_{z^b_0}(i)\text{
and
}(\exists^*z_1^b<\lambda)[\Phi(i,\beta^a,z^a_1,y^b_n,z^b_0,z^b_1)]\right].
\end{equation*}
\end{claim}
\begin{proof}
In (\ref{eqn6}), we can move the quantifier
``$(\forall^*i<\kappa)$'' past the other quantifiers to its left,
so
\begin{multline}
\label{eqn7}
(\forall^*i<\kappa)(\exists^*z^a_1<\lambda)(\exists^*y^b_n<\lambda)\\
(\exists^{*,i}z^b_0<\lambda)(\exists^*
z_1^b<\lambda)[\Phi(i,\beta^a,z^a_1,y^b_n,z^b_0,z^b_1)]
\end{multline}
holds.  The claim will be established if we show that for each
$i<\kappa$ for which
\begin{multline}(\exists^*z^a_1<\lambda)(\exists^*y^b_n<\lambda)\\
(\exists^{*,i}z^b_0<\lambda)(\exists^*
z_1^b<\lambda)[\Phi(i,\beta^a,z^a_1,y^b_n,z^b_0,z^b_1)]
\end{multline}
holds, it is possible to find $v<\lambda_i$ such that
\begin{multline}
(\exists^*z^a_1<\lambda)\biggl{[} g^{\max}_{z^a_1}(i)<v\text{ and }\\
 (\exists^*y^b_n<\lambda)(\exists^*z^b_0<\lambda)
\left[v<f_{z^b_0}(i)\text{ and
}(\exists^*z_1^b<\lambda)[\Phi(i,\beta^a,z^a_1,y^b_n,z^b_0,z^b_1)]\right]\biggr{]}.
\end{multline}
Despite the lengths of the formulas involved, this is not that
hard to accomplish.  Since $\lambda_i<\lambda=\cf(\lambda)$, we
can find $v<\lambda_i$ such that
\begin{multline*}
(\exists^*z^a_1<\lambda)\bigl{[} g^{\max}_{z^a_1}(i)<v\text{ and }\\
(\exists^*y^b_n<\lambda) (\exists^{*,i}z^b_0<\lambda)(\exists^*
z_1^b<\lambda)[\Phi(i,\beta^a,z^a_1,y^b_n,z^b_0,z^b_1)]\bigr{]},
\end{multline*}
and now the result follows from of the definition of
``$\exists^{*,i}z_1^b<\lambda$''.
\end{proof}

Thus there are unboundedly many $z^a_0<\lambda$ for which there is
a function $g\in\prod_{i<\kappa}\lambda_i$ such that for all
sufficiently large $i<\kappa$,
\begin{multline}
\label{eqn10}
(\exists^*z^a_1<\lambda)\biggl{[}g^{\max}_{z^a_1}(i)\leq g(i)\text{ and }\\
(\exists^*y^b_n<\lambda)(\exists^*z^b_0<\lambda)\bigl{[}g(i)<f_{z^b_0}(i)\\
\text{ and
}(\exists^*z_1^b<\lambda)[\Phi(i,z^a_0,z^a_1,y^b_n,z^b_0,z^b_1)]\bigr{]}\biggr{]}.
\end{multline}

Now this is logically equivalent to the statement
\begin{multline}
\label{eqn11}
(\exists^*z^a_1<\lambda)(\exists^*y^b_n<\lambda)(\exists^*z^b_0<\lambda)\\
\bigl{[}g^{\max}_{z^a_1}(i)\leq g(i)<f_{z^b_0}(i)\text{ and
}(\exists^*z_1^b<\lambda)[\Phi(i,z^a_0,z^a_1,y^b_n,z^b_0,z^b_1)]\bigr{]}.
\end{multline}

Suppose we are given a particular $z^a_0<\lambda$ for which a
function $g$ as above can be found, and let us fix $i<\kappa$
``large enough'' so that (\ref{eqn10}) holds.   Also fix  ordinals
$\delta^a<\lambda$ and $\xi^b_n<\lambda$  that serve as suitable
$z^a_1$ and $y^b_n$. Just to be clear, this means that for these
choices we have
\begin{equation*}
(\exists^*z^b_0<\lambda) \bigl{[}g^{\max}_{\delta^a}(i)\leq
g(i)<f_{z^b_0}(i)\text{ and
}(\exists^*z_1^b<\lambda)[\Phi(i,\beta^a,\delta^a,\xi^b_n,z^b_0,z^b_1)]\bigr{]}.
\end{equation*}
Since $\lambda_i<\lambda=\cf(\lambda)$, there must be some value
$w$ satisfying
\begin{equation*}
(\exists^*z^b_0<\lambda) \bigl{[}g(i)<f_{z^b_0}(i)<w\text{ and
}(\exists^*z_1^b<\lambda)[\Phi(i,\beta^a,\delta^a,\xi^b_n,z^b_0,z^b_1)]\bigr{]}.
\end{equation*}
This implies  for our particular $\beta^a$, $g$, $i$, $\delta^a$,
and $\xi^b_n$ that
\begin{multline}
(\forall^* w<\lambda_i)(\exists^*z^b_0<\lambda)\bigl{[}g^{\max}_{\delta^a}(i)\leq g(i)<f_{z^b_1}(i)<w\text{ and }\\
(\exists^*z_1^b<\lambda)[\Phi(i,\beta^a,\delta^a,y^b_n,z^b_0,z^b_1)]\bigr{]}.
\end{multline}
Since $\lambda_i<\lambda=\cf(\lambda)$, the quantifier $(\forall^*
w<\lambda_i)$ can move to the left past the quantifiers
$(\exists^*z^a_1<\lambda)(\exists^*y^b_n<\lambda)$. This tells us
that for our $\beta^a$ and $g$,
\begin{multline}
(\forall^*i<\kappa)(\forall^*w<\lambda_i)(\exists^*z^a_1<\lambda)(\exists^*y^b_n<\lambda)(\exists^*z^b_0<\lambda)\\
\bigl{[}g^{\max}_{z^a_1}(i)\leq g(i)<f_{z^b_0}(i)<w\text{ and }\\
(\exists^*z_1^b<\lambda)[\Phi(i,\beta^a,z^a_1,y^b_n,z^b_0,z^b_1)]\bigr{]}.
\end{multline}
When we put all this together, we end up with the statement
\begin{multline}
(\exists^*z^a_0<\lambda)(\forall^*i<\kappa)(\exists
v<\lambda_i)(\forall^*w<\lambda_i)(\exists^*z^a_1<\lambda)\\
(\exists^*y^b_n<\lambda)(\exists^*z^b_0<\lambda)
\bigl{[}g^{\max}_{z^a_1}(i)\leq v<f_{z^b_0}(i)<w\\
\text{ and
}(\exists^*z_1^b<\lambda)[\Phi(i,\beta^a,z^a_1,y^b_n,z^b_0,z^b_1)]\bigr{]}.
\end{multline}
Since both $\kappa$ and $\lambda_i$ are less than
$\lambda=\cf(\lambda)$, we can move some quantifiers around and
achieve
\begin{multline}
(\forall^*i<\kappa)(\forall^*w<\lambda_i)(\exists^*z^a_0<\lambda)(\exists v<\lambda_i)(\exists^*z^a_1<\lambda)\\
(\exists^*y^b_n<\lambda)(\exists^*z^b_0<\lambda)
\bigl{[}g^{\max}_{z^a_1}(i)\leq v<f_{z^b_0}(i)<w\\
\text{ and
}(\exists^*z_1^b<\lambda)[\Phi(i,\beta^a,z^a_1,y^b_n,z^b_0,z^b_1)]\bigr{]}.
\end{multline}

Thus there is a function $h\in\prod_{i<\kappa}\lambda_i$ such that
\begin{multline}
(\forall^*i<\kappa)(\exists^*z^a_0<\lambda)(\exists v<\lambda_i)(\exists^*z^a_1<\lambda)\\
(\exists^*y^b_n<\lambda)(\exists^*z^b_0<\lambda)
\bigl{[}g^{\max}_{z^a_1}(i)\leq v<f_{z^b_0}(i)<h(i)\\
\text{ and
}(\exists^*z_1^b<\lambda)[\Phi(i,\beta^a,z^a_1,y^b_n,z^b_0,z^b_1)]\bigr{]}.
\end{multline}

After all this work, it is finally time to prove that we can
select objects $\beta^a<\delta^a<\xi^b_n<\beta^b<\delta^b$ that
satisfy all of our requirements.

Clearly, for every unbounded $\Lambda\subs\lambda$,
\begin{equation*}
(\exists
i<\kappa)(\exists^*x\in\Lambda)(h\restr[i,\kappa)<f_x\restr[i,\kappa).
\end{equation*}
Thus we can choose $i^*<\kappa$ such that $h_1^*(i^*)=j_1$ and
$h^*_0(i^*)=n$, and
\begin{multline*}
(\exists^*z^a_0<\lambda)\biggl{[}h\restr[i^*,\kappa)<f^{z^a_0}\restr[i^*,\kappa)\text{ and }(\exists v<\lambda_i) (\exists^*z^a_1<\lambda)(\exists^*y^b_n<\lambda)\\
(\exists^*z^b_0<\lambda)\bigl{[}g^{\max}_{z^a_1}(i^*)\leq v<f_{z^b_0}(i^*)<h(i^*)\text{ and }\\
(\exists^*z^b_1<\lambda)[\Phi(i^*,z^a_0,\dots,z^b_1)]\bigr{]}\biggr{]}.
\end{multline*}

So now we choose $\beta^a$ such that  $h(i^*)<f_{\beta^a}(i^*)$
and for some $\alpha<\lambda_{i^*}$,
\begin{multline*}
(\exists^*z^a_1<\lambda)(\exists^*y^b_n<\lambda)
(\exists^*z^b_0<\lambda)\bigl{[}g^{\max}_{z^a_1}(i^*)\leq \alpha<f_{z^b_0}(i^*)<h(i^*)\text{ and }\\
(\exists^*z^b_1<\lambda)[\Phi(i^*,z^a_0,\dots,z^b_1)]\bigr{]}.
\end{multline*}

Now we choose $\delta^a$, $\xi^b_n$, $\beta^b$, and $\delta^b$
such that
\begin{itemize}
\item $\beta^a<\delta^a<\xi^b_n<\beta^b$
\sk
\item $g^{\max}_{\delta^a}(i^*)\leq\alpha<f_{\beta^b}(i^*)<h(i^*)<f_{\beta^a}(i^*)$
\sk
\item $\Phi(i^*,\beta^a,\delta^a,\xi^b_n,\beta^b,\delta^b)$
\sk
\end{itemize}
It is straightforward to check that these objects satisfy all the
requirements listed in Table 1, so by Claim \ref{final}, we
are done.

\section{Conclusions}

In this final section, we will deduce some conclusions in a few
concrete cases.

\begin{thm}
\label{thm3}
If $\mu$ is a singular cardinal of uncountable cofinality that is not a limit of regular Jonsson cardinals, then $\pr_1(\mu^+,\mu^+,\mu^+,\cf(\mu))$ holds.
\end{thm}

\begin{proof}
The proof of this theorem occurs in two stages---we first show that $\pr_1(\mu^+,\mu^+,\mu,\cf(\mu))$ holds, and then we show that this result can be upgraded to obtain $\pr_1(\mu^+,\mu^+,\mu^+,\cf(\mu)$.

Let $\mu$ be as hypothesized, and let us define $\lambda=\mu^+$ and $\kappa=\cf(\mu)$.

\begin{claim}
$\pr_1(\lambda,\lambda,\mu,\kappa)$ holds.
\end{claim}

\begin{proof}
Let $\langle\kappa_i:i<\kappa\rangle$ be a strictly increasing continuous sequence cofinal in $\mu$. Let $S\subs\{\delta\in [\mu,\lambda):\cf(\delta)=\kappa\}$ be stationary.  Standard  club--guessing results tell us that there is an $S$--club system $\bar{C}$  such that $\id_p(\bar{C},\bar{J})$ is a proper ideal, where $J_\delta$ is the ideal $J^{b[\mu]}_{C_\delta}$ for $\delta\in S$, and furthermore, satisfying $|C_\delta|=\kappa$.  (Note that this last requires that $\kappa=\cf(\mu)$ is uncountable.)

At this point, we have satisfied all of the assumptions of Claim \ref{helpclaim} except possibly for clause (8).  It suffices to show that for each $i<\kappa$, for all sufficiently large regular $\theta<\mu$, Player I has a winning strategy in the game $\gm^{\omega}[\theta,\kappa_i,1]$.  Since $\mu$ is not a limit of regular Jonsson cardinals, it follows that for all sufficiently large regular $\theta<\mu$, Player I has a winning strategy in $\gm^\omega[\theta,\theta,1]$.  This implies, by Lemma \ref{1.4} (1), that for all sufficiently large regular $\theta$, Player I has a winning strategy in $\gm^\omega[\theta,\kappa_i,1]$, and so clause (8) of Claim \ref{helpclaim} is satisfied.
\end{proof}

To finish the proof of Theorem \ref{thm3}, it remains to show that we can increase the number of colors from $\mu$ to $\lambda=\mu^+$ --- we need $\pr_1(\lambda,\lambda,\lambda,\kappa)$ instead of $\pr_1(\lambda,\lambda,\mu,\kappa)$.

\begin{lem}
There is a coloring $c_1:[\lambda]^2\rightarrow\lambda$ such that whenever we are given
\begin{itemize}
\item $\theta<\kappa$,
\sk
\item $\langle t_\alpha:\alpha<\lambda\rangle$ a sequence of pairwise disjoint elements of $[\lambda]^\theta$,
\sk
\item $\zeta_\alpha\in t_\alpha$ for $\alpha<\lambda$, and
\sk
\item $\Upsilon<\lambda$,
\sk
\end{itemize}
we can find $\alpha<\beta$ such that $t_\alpha\subs\min(t_\beta)$ and
\begin{equation}
(\forall\zeta\in t_\alpha)[c_1(\zeta,\zeta_\beta)=\Upsilon].
\end{equation}
\end{lem}
\begin{proof}
Let $c:[\lambda]^2\rightarrow\mu$ be a coloring that witnesses $\pr_1(\lambda,\lambda,\mu,\kappa)$.  For each $\alpha<\lambda$, let $g_\alpha$ be a one--to--one function from $\alpha$ into $\mu$.   We define
\begin{equation}
c_1(\alpha,\beta)=g^{-1}_\beta(c(\alpha,\beta)).
\end{equation}

Suppose now that we are given objects $\theta$, $\langle t_\alpha:\alpha<\lambda\rangle$, $\langle \zeta_\alpha:\alpha<\lambda\rangle$, and $\Upsilon$ as in the statement of the lemma.  Clearly we may assume that $\min(t_\alpha)>\alpha$.

For $i<\mu$, we define $X_i:=\{\alpha\in[\gamma,\lambda):g_{\zeta_\alpha}(\Upsilon)=i\}$.  Since $\lambda$ is a regular cardinal, it is clear that there is $i^*<\mu$ for which $|X_{i^*}|=\lambda$.    Since $c$ exemplifies $\pr_1(\lambda,\lambda,\mu,\kappa)$, for some $\alpha<\beta$ in $X_{i^*}$ we have $t_\alpha\subs \min(t_\beta)$ and
\begin{equation}
(\forall\zeta\in t_\alpha)[c(\zeta,\zeta_\beta)=i^*].
\end{equation}
By definition, this means
\begin{equation}
(\forall\zeta\in t_\alpha)[c_1(\zeta,\zeta_\beta)=g^{-1}(c(\alpha,\beta))=g^{-1}(i^*)=\Upsilon],
\end{equation}
hence $\alpha$ and $\beta$ are as required.

\end{proof}

To continue the proof of Theorem \ref{thm3}, we define a coloring $c_2:[\lambda]^2\rightarrow\lambda$ by
\begin{equation}
c_2(\alpha,\beta)=c_1(\alpha,\nu(\alpha,\beta)),
\end{equation}
where $\nu(\alpha,\beta)$ is as in the proof of Theorem \ref{ash}.

It remains to check that $c_2$ witnesses $\pr_1(\lambda,\lambda,\lambda,\kappa)$.  Toward this end, suppose we are given $\theta<\kappa$, $\langle t_\alpha:\alpha<\lambda\rangle$ a sequence of pairwise disjoint members of $[\lambda]^\theta$, and $\Upsilon<\lambda$.  We need to find $\delta^a$ and $\delta^b$ less than $\lambda$ such that
\begin{equation}
\epsilon^a\in t_{\delta^a}\wedge\epsilon^b\in t_{\delta^b}\Longrightarrow c_2(\epsilon^a,\epsilon^b)=\Upsilon.
\end{equation}

\begin{lem}
There is a stationary set of $\gamma_1<\lambda$ such that for some $\gamma_0<\gamma_1$ and $\beta\in[\gamma_1,\lambda)$, if $\gamma_0\leq\alpha<\gamma_1$, then the function $\nu$ is constant on $t_\alpha\times t_\beta$.
\end{lem}
\begin{proof}
Let $E$ be an arbitrary closed unbounded subset of $\lambda$, and let $W$ be the set of ordinals $<\lambda$ satisfying the properties of $\gamma_1$.  In the proof of Theorem \ref{ash}, without loss of generality we can have $E\in M_0$.   This means that the ordinal $\beta^*$ found in the course of that proof will be in $E$, so we finish by observing that $\beta^*\in W$.
\end{proof}   

An application of Fodor's Lemma gives us a single ordinal $\gamma_0$ and a stationary $W'\subs W$ such that for all $\gamma\in W'$, there is a $\beta_\gamma\in[\gamma,\lambda)$ such that for all $\alpha\in [\gamma_0,\gamma)$, $\nu\restr (t_\alpha\times t_\beta)$ is constant.

Using properties of the coloring $c_1$, we can find $\alpha$ and $\gamma$ such that
\begin{itemize}
\sk
\item $\gamma_0\leq\alpha<\lambda$
\sk
\item $\gamma\in W'\setminus (\sup(t_\alpha)+1)$, and
\sk
\item $\zeta\in t_\alpha\Longrightarrow c_1(\zeta,\gamma)=\Upsilon$.
\sk
\end{itemize}
Now given $\epsilon^a\in t_\alpha$ and $\epsilon^b\in t_{\beta_\gamma}$, we find
\begin{equation}
c_2(\epsilon^a,\epsilon^b)=c_1(\epsilon^a,\gamma)=\Upsilon,
\end{equation}
and therefore $c_2$ exemplifies $\pr(\lambda,\lambda,\lambda,\kappa)$.
\end{proof}

Theorem \ref{thm3} strengthens results in  \cite{Sh:413} as clearly $\pr_1(\mu^+,\mu^+,\mu^+,\cf(\mu))$ implies that $\mu^+$ has a Jonsson algebra (i.e., $\mu^+$ is not a Jonsson cardinal).  The question of whether the successor of a singular cardinal can be a Jonsson cardinal is a well--known open question.

We note that many of the results from Section 2 of \cite{Sh:413} dealing with the existence of winning strategies for Player I in $\gm^\omega[\lambda,\mu,\gamma]$ can be combined with Theorem \ref{ash} to give new results. For example, we have the following result from \cite{Sh:413}.

\begin{prop}  
\label{413prop}  
If $\tau\leq 2^\kappa$ but $(\forall \theta<\kappa)[2^\theta<\tau]$, then Player I has a winning strategy in the game $\gm^\omega(\tau,\kappa,\kappa^+)$.
\end{prop}
\begin{proof}
See Claim 2.3(1) and Claim 2.4(1) of \cite{Sh:413}.
\end{proof}

Armed with this, the following claim is straightforward.

\begin{claim}
Let $\mu$ be a singular cardinal of uncountable cofinality.  Further assume that $\chi$ is a cardinal such that
$2^{<\chi}\leq\mu<2^{\chi}$.  Then $\pr_1(\mu^+,\mu^+,\chi,\cf(\mu))$ holds.
\end{claim}
\begin{proof}
If $2^{<\chi}<\mu$, then Claims 2.3(1) and 2.4(1) of \cite{Sh:413} imply that for every sufficiently large $\theta<\mu$, Player I has a winning strategy in the game $\gm^\omega(\theta,\chi,\chi^+)$.

If $\mu=2^{<\chi}$, then $\cf(\mu)=\cf(\chi)$. Let $\langle\kappa_i:i<\cf(\mu)\rangle$ be a strictly increasing continuous sequence of cardinals cofinal in $\chi$.  Given $i<\cf(\mu)$, we claim that for all sufficiently large regular $\tau<\mu$, Player I has a winning strategy in $\gm^\omega(\tau,\kappa_i,\chi)$.  Once we have established this, $\pr_1(\mu^+,\mu^+,\chi,\cf(\mu))$ follows by Theorem \ref{ash}.

Given $\tau=\cf(\tau)$ satisfying $2^{\kappa_i}<\tau<\mu$, let $\eta$ be the least cardinal such that $\tau\leq 2^\eta$.  Clearly $\kappa_i<\eta<\chi$.  By Proposition \ref{413prop}, Player I wins the game $\gm^\omega(\tau,\eta,\eta^+)$.  This implies (since $\eta^+<\chi$ and $\kappa_i<\eta$) that Player I wins the game $\gm^\omega(\tau,\kappa_i,\chi)$ as required.
\end{proof}

We can also use Claim \ref{1.4b} to prove similar results. For example we have the following.

\begin{claim}
Let $\mu$ be a singular cardinal of uncountable cofinality.  Further assume that $\chi<\mu$ satisfies $2^\chi<\mu<\beth_{(2^\chi)^+}(\chi)$.  Then $\pr_1(\mu^+,\mu^+,\chi,\cf(\mu))$ holds.
\end{claim}

\begin{proof}
Again, the main point is that for all sufficiently large regular $\theta<\mu$, Player I has a winning strategy in the game $\gm^\omega[\theta,\chi,(2^\chi)^+]$. This follows immediately from Claim \ref{1.4b}.  Since $(2^\chi)^+<\mu$, Theorem \ref{ash} is applicable.
\end{proof}

In a sequel to this paper, we will address the situation where $\lambda$ is the successor of a singular cardinal of countable cofinality.  Similar results hold, but the combinatorics involved are trickier.

\def\germ{\frak} \def\scr{\cal}
  \ifx\documentclass\undefinedcs\def\rm{\fam0\tenrm}\fi%f**k-amstex!
  \def\defaultdefine#1#2{\expandafter\ifx\csname#1\endcsname\relax
  \expandafter\def\csname#1\endcsname{#2}\fi} \defaultdefine{Bbb}{\bf}
  \defaultdefine{frak}{\bf} \defaultdefine{mathbb}{\bf}
  \defaultdefine{mathcal}{\cal}
  \defaultdefine{beth}{BETH}\defaultdefine{cal}{\bf} \def\bbfI{{\Bbb I}}
  \def\mbox{\hbox} \def\text{\hbox} \def\om{\omega} \def\Cal#1{{\bf #1}}
  \def\pcf{pcf} \defaultdefine{cf}{cf} \defaultdefine{reals}{{\Bbb R}}
  \defaultdefine{real}{{\Bbb R}} \def\restriction{{|}} \def\club{CLUB}
  \def\w{\omega} \def\exist{\exists} \def\se{{\germ se}} \def\bb{{\bf b}}
  \def\equivalence{\equiv} \let\lt< \let\gt> \def\cite#1{[#1]}
\providecommand{\bysame}{\leavevmode\hbox
to3em{\hrulefill}\thinspace}
\providecommand{\MR}{\relax\ifhmode\unskip\space\fi MR }
% \MRhref is called by the amsart/book/proc definition of \MR.
\providecommand{\MRhref}[2]{%
  \href{http://www.ams.org/mathscinet-getitem?mr=#1}{#2}
} \providecommand{\href}[2]{#2}

\end{document}